
\documentstyle{amsppt}
\topmatter

\parskip=0.5\baselineskip
\baselineskip=1.1\baselineskip
\parindent=0pt
\loadmsbm
\UseAMSsymbols
\hoffset=0.75truein
\voffset=0.5truein

\font\twrm=cmr8
\def\doct{\delta_{oct}}
\def\pt{\hbox{\it pt}}

\def\sol{\operatorname{sol}}
\def\dih{\operatorname{dih}}
\def\vor{\operatorname{vor}}
\def\quo{\operatorname{quo}}
\def\PM{\operatorname{PM}}
\def\octavor{\operatorname{octavor}}

\def\rad{\operatorname{rad}}

\def\refy{\relax}
\def\tlp{\tau_{\hbox{\twrm LP}}}  
\def\squander{(4\pi\zeta-8)\,\pt}

\def\diag|#1|#2|{\vbox to #1in {\vskip.3in\centerline{\tt Diagram #2}\vss} }
\def\v{\hskip -3.5pt }
\def\gram|#1|#2|#3|{
        {
        \smallskip
        \hbox to \hsize
        {\hfill
        \vrule \vbox{ \hrule \vskip 6pt \centerline{\it Diagram #2}
         \vskip #1in %
             \includegraphics{#3}\hrule }
        \v\vrule\hfill
        }
\smallskip}}

\title Sphere Packings III\endtitle
\author Thomas C. Hales\endauthor
\endtopmatter

\abstract{Abstract}
An earlier paper describes a five-step program to prove the Kepler
Conjecture.  This paper carries out the third step of the program.
\endabstract

\document
\footnote""{\line{\hfill\it version - 7/18/98, revised 1/19/02}}
\footnote""{Research partially funded by the NSF}

\head 1. Introduction and Review\endhead

This paper is a continuation of the first two parts of this
series (\cite{I},\cite{II}). It relies on the formulation of the
Kepler conjecture in \cite{F}. The terminology and notation of this
paper are consistent with these earlier papers, and we refer
to results from them by prefixing the relevant
section numbers with I, II, or F. Around each vertex
is a modification of the Voronoi cell, called the $V$-cell and
a collection of quarters and quasi-regular tetrahedra.
These objects constitute the {\it decomposition star\/} at the vertex.
 A decomposition star may be decomposed into {\it
standard clusters}.
By definition, a
standard cluster is the part
of the given decomposition star
that lies over a given {\it standard region\/}
on the unit sphere.

A real number, called the {\it score},
 is attached to each cluster.
Each star
receives a score by summing the scores $\sigma(R)$ for the clusters
$R$ in the
star.   The scores are measured in multiples of a {\it point\/}
(\pt), where
$\pt \approx 0.055$.  If every star scores at most $8\,\pt$,
then the Kepler conjecture follows.

The steps of the Kepler conjecture, as outlined in Part I, are

{

\def\ha{ \hangindent=20pt \hangafter=1\relax }
1. A proof that if all standard
regions are triangular, the total score
is less than $8\,\pt$

\ha
2.  A proof that the standard regions
with more than
three sides
score at most $0\,\pt$

\ha
3. A proof that if all of the
standard regions are triangles or quadrilaterals,
then the total score is less than $8\,\pt$ (excluding the
case of pentagonal prisms)

\ha
4.  A proof that if some
standard region has
more than four sides, then the
star scores less than $8\,\pt$

\ha
5.  A proof that pentagonal prisms score less than $8\,\pt$

}

This paper completes the third of these steps.

The standard regions of the decomposition
stars in the face-centered cubic and the hexagonal-close packings
are regular
triangles and quadrilaterals.  These stars score exactly
$8\,\pt$.  The local optimality results for hexagonal-close packings
and face-centered cubic packings have been established in \cite{II}.
If the planar map is that of a close packing, there are eight quasi-regular
tetrahedra.  To score $8\,\pt$, the quasi-regular tetrahedra must
be regular of edge length 2 \cite{I.9.1}.  Such decomposition stars
are precisely those of the close packings.  Conjecturally,
all the decomposition stars have scores strictly less than $8\,\pt$,
so that various approximations may be introduced to prove the desired
bounds.

The standard regions of pentagonal prisms are
triangles and quadrilaterals (10 triangles and 5
quadrilaterals).  The pentagonal prisms are the subject of
the fifth step of the program.

\proclaim{Theorem 1}
 Let $D^*$ be a decomposition star whose combinatorial
structure is not a pentagonal prism.
Suppose that
each standard region of $D^*$ is a triangle or quadrilateral.
Then the score
of $D^*$ is at most $8\,\pt$.  Equality is attained exactly
when the decomposition of the unit sphere into standard regions
coincides with that of a decomposition star in
the face-centered cubic or hexagonal-close packing.
\endproclaim

The proof of Theorem 1 relies on many computer calculations.
We make a list of combinatorial properties that a decomposition
star must have for it to have a possibility of scoring more
than $8\,\pt$.  We then make a computer
search to find all decompositions of the unit sphere into
triangles and quadrilaterals that satisfy all the properties
on the list.  The computer search produces an explicit
list containing nearly 2000
combinatorial types.

For each of these combinatorial arrangements
of triangles and quadrilaterals,
we have a nonlinear optimization problem:  maximize the
score over the space of all decomposition stars $D^*$ with
the given combinatorial arrangement.
It is not necessary to solve this optimization problem. It is
sufficient to establish an upper bound of $8\,\pt$.  To do this,
we define a {\it linear relaxation\/} of the original problem,
that is, a linear programming
problem whose solution strictly dominates the
global maximum of the original nonlinear problem.
This gives an upper
bound on the linear problem,  which is usually less than
$8\,\pt$.

In some cases, this procedure leads to a bound
greater than $8\,\pt$, and further analysis will be required.

One advantage of our method of linear relaxation
is that the verification of
the bounds is particularly simple.  If the linear relaxation
asks to maximize $c\cdot x$ subject to the system of linear inequalities
$Ax\le b$ and  $x\ge0$, then duality theory produces a vector $z$
with nonnegative entries such that $c\le zA$.
To verify the bound of $8\,\pt$, it is enough to
check that  $c\le zA$, $z\ge0$, and  $z\cdot b<8\,\pt$ (because
then $c\cdot x \le z A x \le z\cdot b < 8\,\pt$).

\medskip
Many of the linear inequalities that were used in the
linearly relaxed optimization problems were obtained as
follows.  We use numerical methods to find a convex
polygon containing the set of ordered pairs
$$(\dih(S), \sigma(S))$$
as $S$ ranges over quasi-regular tetrahedra.
  The edges of the polygon give linear inequalities
relating the dihedral angles to the score $\sigma$.
More generally, additional inequalities are obtained by considering
polygons that contain the ordered pairs
$$(\dih(S),\sigma(S)-\lambda\sol(S)),$$
for appropriate constants $\lambda$.   Similar inequalities
are obtained for quadrilateral regions.

We must remember that ultimately
we are dealing with a nonlinear optimization problem
that is larger, by a considerable order,
than what is conventionally thought to
be solvable by exact methods.  The domain of our
optimization problem has many components,
and the dimension varies from component
to component. Even the magnitude of the problem is poorly understood.  The
best-known bound on the
dimension of the components is about 155 dimensions.  Components
of interest
frequently have more than 35 dimensions.
 The decomposition into standard
regions gives the problem an aspect of separability.
Nevertheless, certain complications will have
to be tolerated.

The biggest weakness of this method is that the output from
the computer search for the combinatorial arrangements of
triangles and quadrilaterals is not easily checked for errors.
The algorithm is described in some detail in Section 9, but the
only assurance that no cases have been skipped comes through a
careful reading of the computer code.
It would be advantageous to have a more transparent proof of the
results of this section.

This paper is supplemented by an
appendix giving further details of the combinatorial arrangements.
Additional
details about these calculations, including the full source
code for all of the computer verifications of this paper, can be
found in \cite{H2}.

The outline of this paper was developed at the University of Chicago
during the summer of 1994.  I would like to
thank to P. Sally for making computer resources and
other facilities available to me during my stay in Chicago,
and for his encouragement with this work.
I would also like to give special thanks to S. Ferguson for many helpful
discussions concerning this topic.   His investigations have led
to a number of improvements in the results presented here.

\head 2. Geometric Considerations\endhead

\bigskip
\subhead 2.1\endsubhead
We will call a standard cluster over a quadrilateral
region a {\it quad cluster}.
The four vertices of the quad cluster whose projections to the
unit sphere mark the
extreme points of the quadrilateral region will be called
the {\it corners\/} of the cluster.
We call the four angles of the
standard region associated with the quad cluster its {\it dihedral\/}
angles.

The rules defining the score
have undergone a long series of revisions
over the last several years.  The formulation used in this
paper is described in \cite{F}.

\proclaim {Lemma 2.2}
A quadrilateral region does not enclose any
vertices of height $2.51$ or less.
\endproclaim

\demo{Proof}
Let $v_1,\ldots,v_4$ be the corners of the quad cluster, and
let $v$ be an enclosed vertex of height at most $2.51$.
We cannot have $|v_i-v|\le2.51$ for two different vertices
$v_i$,  because two
such inequalities would partition the region into
 two separate
standard regions instead of a single quadrilateral region.
We apply I.4.3 to simplify the quad cluster.  (Lemma I.4.3 assumes
the existence of another enclosed vertex $v'$, but it can be
omitted both from the statement of the Lemma and from the proofs without
affecting matters.)  Then I.4.3 allows us to assume
$$|v_i-v_{i+1}|=2.51,\quad |v_i|=2, \quad |v|=2.51,$$
for $i=1,\ldots,4$.
Reindexing and perturbing $v$ as necessary,
we may assume that $2\le |v_1-v|\le2.51$
and $|v_i-v|\ge2.51$, for $i=2,3,4$.
Moving $v$, we may assume it reaches the minimal distance to two adjacent
corners ($2$ for $v_1$ or $2.51$ for $v_i$, $i>1$).  Keeping $v$
fixed at this minimal distance, perturb the quad cluster along its
remaining degree of freedom until $v$ attains its minimal distance
to three of the corners.  This is a rigid figure.  There are four
possibilities depending on which three corners are chosen.
Pick coordinates to show that the distance from $v$ to the remaining
vertex violates its inequality.\qed
\enddemo

\head 3. Functions Related to the Score\endhead

Set $\zeta^{-1}:=\sol(S(2,2,2,2,2,2))=2\arctan(\sqrt{2}/5)$.
The constant $\zeta$ is related to the other fundamental
constants by the relations $\pt= 2/\zeta-\pi/3$ and
$\doct=(\pi-2/\zeta)/\sqrt{8}$.  Rogers's bound
is $\sqrt{2}/\zeta\approx 0.7796$.

We consider the functions
$\sigma_\lambda(R):=\sigma(R)-\lambda \zeta\sol(R)\,\pt$, for
$\lambda=0$, $1$, or $3.2$, where
$R$ is a standard cluster.
The constant $3.2$ was determined experimentally.
We will see that $\sigma_1$ has a simple interpretation.
 We write $\tau(R) = -\sigma_1(R)$.
If $D^*$ is a decomposition star with standard clusters $\{R\}$, set
$\tau(D^*) = \sum_{R}\tau(R)$.
\smallskip

\proclaim{Lemma 3.1}
$\tau(R)\ge 0$, for all standard clusters $R$.
\endproclaim

\demo{Proof}  If $R$ is not a quasi-regular tetrahedron, or
if it is but $\rad(R)\ge1.41$, then
$\sigma(R)\le0$ and $\sol(R)> 0$, so that the result is immediate
(see I.9.17).
Assume that $R$ is a
quasi-regular tetrahedron and $\rad(R)\le 1.41$.
  The result follows from
Calculation 10.3.6, which asserts that $\Gamma(R)\le \sol(R)\zeta\pt$.
(Equality is attained only for the regular tetrahedron of edge $2$.)
  \qed
\enddemo

\proclaim{Lemma 3.2}
$$\sigma(D^*) = {4\pi \zeta\,\pt} - \tau(D^*).$$
\endproclaim

\demo{Proof}
Let $R_1,\ldots,R_k$ be the standard clusters in $D^*$. Then
$$\sigma(D^*) = \sum\sigma(R_i) + (4\pi-\sum\sol(R_i))\zeta\,\pt = 4\pi \zeta\,\pt - \sum\tau(R_i).$$
\qed
\enddemo

Since $4\pi \zeta< 22.8$, we find as an immediate corollary that
if there are standard clusters satisfying $\tau(R_1)+\cdots+\tau(R_k)\ge14.8\,\pt$,
then the score of the star is less than $8\,\pt$.

The function $\tau(R)$ gives the amount {\it squandered\/} by
a particular standard cluster $R$.  If nothing is squandered,
then $\tau(R_i)=0$ for every standard cluster, and the upper bound
is $4\pi \zeta\,\pt\approx 22.8\,\pt$.
This is Rogers's bound
on density.  It is the unattainable bound that would
be obtained by
packing regular tetrahedra around a common vertex with no distortion and
no gaps.
(More precisely, in the terminology of \cite{H1}, the score
$s_0=4\pi \zeta\,\pt$
corresponds to the {\it effective density\/}
 $16\pi\doct/(16\pi- 3 s_0)  =\sqrt{2}/\zeta \approx 0.7796$,
which is Rogers's bound.)  Every positive lower bound
on $\tau(R_i)$ translates into an improvement on Rogers's bound.
 To say that a decomposition star scores at most $8\,\pt$
is to say that at least
$(4\pi \zeta-8)\pt\approx 14.8\,\pt$ are squandered.

\head 4. Some Linear Constraints\endhead

This section gives some linear inequalities between
$\sigma(R)-\lambda \sol(R)\zeta\pt$ and
$\dih(R)$.

\proclaim {Proposition 4.1}  Let $R$ be a quad cluster.
  Let $\sigma(R)$ denote its score,
let $\dih(R)$ be one of the four dihedral angles of $R$,
 and let $\sol(R)$ be the solid angle of the
standard region of $R$.  The following inequalities hold among
$\dih(R)$, $\sol(S)$,  and $\sigma(R)$:
\endproclaim

{\parskip=0pt
\hbox{}

$1$: \quad $\sigma(R)< -5.7906 + 4.56766 \dih(R)$,

$2$: \quad $\sigma(R)< -2.0749 + 1.5094 \dih(R)$,

$3$: \quad $\sigma(R)< -0.8341 + 0.5301 \dih(R)$,

$4$: \quad $\sigma(R) < -0.6284 + 0.3878\dih(R)$,

$5$: \quad $\sigma(R) < 0.4124 - 0.1897 \dih(R)$,

$6$: \quad $\sigma(R) < 1.5707- 0.5905\dih(R)$,

$7$: \quad $\sigma(R) < 0.41717 - 0.3\sol(R)$,

$8$: \quad $\sigma_1(R) < -5.81446 + 4.49461 \dih(R)$,

$9$: \quad $\sigma_1(R) < -2.955 + 2.1406 \dih(R)$,

$10$: \quad $\sigma_1(R) < -0.6438 + 0.316 \dih(R)$,

$11$: \quad $\sigma_1(R) < -0.1317$,

$12$: \quad $\sigma_1(R) < 0.3825 - 0.2365 \dih(R)$,

$13$: \quad $\sigma_1(R) < 1.071 - 0.4747 \dih(R)$,

$14$: \quad $\sigma_{3.2}(R) < -5.77942 + 4.25863\dih(R)$,

$15$: \quad $\sigma_{3.2}(R) < -4.893 + 3.5294 \dih(R)$,

$16$: \quad $\sigma_{3.2}(R) < -0.4126$,

$17$: \quad $\sigma_{3.2}(R) < 0.33 - 0.316 \dih(R)$,

$18$: \quad $\sigma(R) < -0.419351 \sol(R) -5.350181+ 4.611391\dih(R)$,

$19$: \quad
    $\sigma(R) < -0.419351 \sol(R) -1.66174 + 1.582508\dih(R)$,

$20$: \quad
    $\sigma(R) < -0.419351 \sol(R) +0.0895+ 0.342747\dih(R)$,

$21$: \quad
    $\sigma(R) < -0.419351 \sol(R) +3.36909 - 0.974137\dih(R)$.

}

\proclaim{Proposition 4.2}
 Let $R$ be a quad cluster.
  Let $\dih_1(R)$ and $\dih_2(R)$
be two adjacent dihedral angles of $R$.
  Set $d(R) = \dih_1(R)+\dih_2(R)$.  The following
inequalities hold between $d(R)$ and $\sigma(R)$.
\endproclaim

{\parskip=0pt

\hbox{}

$1$:  \quad $\sigma(R) < -9.494 + 3.0508\, d(R)$,

$2$:  \quad $\sigma(R) < -1.0472 + 0.27605\, d(R)$,


$3$:  \quad $\sigma(R) < 3.5926 - 0.844 \, d(R)$,

}

\proclaim {Proposition \refy{4.3}}
\endproclaim

{\parskip=0pt
\hbox{}

$1$: \quad $1.153< \dih(R)$,

$2$: \quad $\dih(R)< 3.247$.

}

\demo{Proof}  Proposition 4.3 follows from
interval arithmetic calculations
based on the methods of \cite{I}.
The lower bound of 1.153 holds in fact for the dihedral angles
of any standard cluster other than quasi-regular tetrahedra.
\qed
\enddemo

All of these inequalities have been proved by
interval arithmetic methods
by computer.
An appendix gives
a general description of the cases involved in the verification.
  Further details are available \cite{H2}.

\head 5. Types of Vertices\endhead

The combinatorial structure of a decomposition star is conveniently
described as  a {\it planar map}.
A planar graph is a graph that can be embedded into the
plane or sphere.  A planar map is a planar graph with
additional combinatorial structure that encodes a particular embedding
of the graph \cite{T}.  All our planar maps will be unoriented:
we do not distinguish between a planar map and its reflection.
Associated with a planar map are faces, (combinatorial) angles
between adjacent edges, and so forth.  Associated with each planar
map $L$ is a planar graph $G(L)$, obtained by forgetting the additional
combinatorial structure.  Each planar map has a dual $L^*$, obtained by
interchanging faces and vertices.
The faces of a planar map are in natural
bijection with the vertices of $L^*$.  We say that a face
is an $n$-gon if the corresponding vertex in the dual $L^*$
has degree $n$.  The {\it boundary\/} of a face is
an $n$-circuit in $G(L)$. The edges of the boundary are in natural bijection
with the edges in $L^*$ that are joined to the
dual vertex of the face in $L^*$.

Associated with each decomposition star is a standard decomposition
of the unit sphere, as described in Part I.  We form
a planar map $L$ by associating with each standard region a face of $L$
and with each edge of a standard region an edge of $L$.
This paper is concerned with the special case of the Kepler
conjecture in which every face of $L$ is a triangle or quadrilateral.

We say that a vertex $v$ of $L$ has {\it type\/} $(p,q)$ if there
are exactly $p$ triangular faces and $q$ quadrilateral faces
that meet at $v$. We write $(p_v,q_v)$ for the type of $v$.
The type of a vertex of the decomposition star is the type
of the corresponding vertex in the planar map.

We use the following strategy in the proof of step 3 of the
Kepler conjecture.  The linear inequalities that were stated
in Section 4 will be combined to give a bound on the score of
the standard clusters around a given vertex of a given type.
This bound will depend only on the type of the vertex.
The bound comes as the solution to the linear programming
problem of optimizing
the sum of scores, subject to the linear constraints of Section 4
and to the constraint that the dihedral angles around the vertex
sum to $2\pi$.  Similarly, we obtain a lower bound on what is
squandered around each vertex.

This gives certain obvious
constraints on decomposition stars.  For example, if more than
$\squander$ are squandered at a vertex of a given type,
then that type of vertex cannot be part of a decomposition star
scoring more than $8\,\pt$.  These relations between
scores and vertex types
will allow us to reduce the feasible planar maps to an
explicit finite list.
For each of the planar maps on this list, we calculate
a second, more refined linear programming bound on the
score.
Often, the refined linear programming bound is less than $8\,\pt$.

This section derives the bounds on the scores of the
clusters around a given vertex as a function of the
type of the vertex.  Define constants
$\tlp(p,q)/\pt$ by Table 5.1.  The entries marked with an asterisk will
not be needed.

\bigskip

\def\pt{\hbox{\it pt}}

$$
\vbox{\offinterlineskip
\hrule
\halign{&\vrule#&\strut\ \hfil#\hfil\ \cr   
height 7pt&\omit&&\omit&&\omit&&\omit&&\omit&&\omit&&\omit&\cr
&\hfil $\tlp(p,q)/\pt$\hfil
        &&\hfil $q=0$\hfil
        &&\hfil1\hfil
        &&\hfil2\hfil
        &&\hfil3\hfil
        &&\hfil4\hfil
        &&\hfil5\hfil&
\cr
height 7pt&\omit&&\omit&&\omit&&\omit&&\omit&&\omit&&\omit&\cr
\noalign{\hrule}
height7pt&\omit&&\omit&&\omit&&\omit&&\omit&&\omit&&\omit&\cr
&$p=0$&&    *&& *&& 15.18&& 7.135&& 10.6497&& 22.27&\cr
&1&&    *&& *&&  6.95&& 7.135&&17.62  && 32.3&\cr
&2&&    *&& 8.5&&4.756&&12.9814&&*&&*&\cr
&3&&    *&& 3.6426&&8.334&&20.9&&*&&*&\cr
&4&&4.1396&&3.7812&&16.11&&*&&*&&*&\cr
&5&&0.55&&11.22&&*&&*&&*&&*&\cr
&6&&6.339&&*&&*&&*&&*&&*&\cr
&7&&14.76&&*&&*&&*&&*&&*&\cr
height7pt&\omit&&\omit&&\omit&&\omit&&\omit&&\omit&&\omit&\cr}
\hrule
}\tag 5.1
$$

\bigskip

\proclaim{Proposition 5.2}
Let $S_1,\ldots,S_p$ and $R_1,\ldots,R_q$ be
the tetrahedra and quad clusters around a vertex of type $(p,q)$.
Consider the constants
of Table 5.1.  We have
$$\align
&\sum^p\tau(S_i) + \sum^q\tau(R_i) \ge \tlp(p,q),\\
\endalign
$$
\endproclaim

\demo{Proof}
Set
$$(d_i^0,t_i^0)=(\dih(S_i),\tau(S_i)),\qquad
(d_i^1,t_i^1)=(\dih(R_i),\tau(R_i)).$$  Then
$\sum^p\tau(S_i)+\sum^q\tau(R_i)$ is at least the minimum
of $\sum^p t_i^0+\sum^q t_i^1$ subject to
$\sum^p d_i^0+\sum^q d_i^1 = 2\pi$ and to the system
of linear inequalities of Section 10 (Group 3) and
Proposition 4.1 (obtained
by replacing $-\sigma_1$ and dihedral angles by $t_i^j$ and $d_i^j$).
The constant $\tlp(p,q)$ was chosen to be slightly larger
than the actual minimum of this linear programming problem.

The entry $\tlp(5,0)$ is based on Lemma
5.3, $k=1$.  \qed
\enddemo

\proclaim{Lemma 5.3}
Let $v_1,\ldots, v_k$, for some
$k\le 4$, be distinct vertices of a decomposition
star of type $(5,0)$.  Let $S_1,\ldots, S_r$ be quasi-regular
tetrahedra around the edges $(0,v_i)$, for $i\le k$.
Then
$$\sum_{i=1}^r \tau(S_i)> 0.55k\,\pt,$$
and
$$\sum_{i=1}^r \sigma(S_i) < r\,\pt - 0.48k\,\pt.$$
\endproclaim

\demo{Proof}
We have $\tau(S)\ge 0$, for any quasi-regular
tetrahedron $S$.  We refer to the edges $y_4,y_5,y_6$ of a simplex
$S(y_1,\ldots,y_6)$ as its top edges. Set $\xi=2.1773$.

We claim (Claim 1) that if $S_1,\ldots,S_5$ are quasi-regular tetrahedra around
an edge $(0,v)$ and if $S_1=S(y_1,\ldots,y_6)$, where $y_5\ge\xi$
is the length of a top edge $e$ on $S_1$ shared with $S_2$, then
$\sum_1^5\tau(S_i) > 3(0.55)\,\pt$.  This claim follows from Inequalities
10.5.1 and 10.5.2 if some other top edge in this group
of quasi-regular tetrahedra has length greater than $\xi$.
Assuming all the top edges other than $e$ have length at most
$\xi$, the estimate follows from $\sum_1^5\dih(S_i)=2\pi$ and
Inequalities 10.5.3, 10.5.4.

Now let $S_1,\ldots,S_8$ be the eight quasi-regular tetrahedra
around two edges $(0,v_1)$, $(0,v_2)$ of type $(5,0)$.
Let $S_1$ and $S_2$ be the simplices along the face $(0,v_1,v_2)$.
Suppose
that the top edge $(v_1,v_2)$ has length at least $\xi$.
We claim (Claim 2) that $\sum_1^8\tau(S_i)> 4(0.55)\,\pt$.  If there is a top
edge of length at least $\xi$ that does not lie on $S_1$ or
$S_2$,
then this claim reduces to
Inequality 10.5.1 and Claim 1.
If any of the top edges of $S_1$ or $S_2$ other than $(v_1,v_2)$ has
length at least $\xi$, then the claim follows from
Inequalities 10.5.1 and 10.5.2.  We assume all top edges other
than $(v_1,v_2)$ have length at most $\xi$.  The claim now
follows from Inequalities 10.5.3 and 10.5.5, since the dihedral
angles around each vertex sum to $2\pi$.

We prove the bounds for $\tau$.  The proof for $\sigma$ is entirely
similar, but uses the constant $\xi=2.177303$ and the
Inequalities 10.5.8--10.5.14 rather than 10.5.1--10.5.7.
Claims analogous to Claims 1 and 2 hold for the $\sigma$ bound
by Inequalities 10.5.8--10.5.12.

Consider $\tau$ for $k=1$.  If a top edge has length at least $\xi$, this
is Inequality 10.5.1.  If all top edges have length less than $\xi$,
this is Inequality 10.5.3, since dihedral angles sum to $2\pi$.

We say that a top edge lies around a vertex $v$ if it is an
edge of a quasi-regular tetrahedron with vertex $v$.
We do not require $v$ to be the endpoint of the edge.

Take $k=2$.
If there is an edge of length at least $\xi$ that
lies around only one of $v_1$ and $v_2$, then Inequality 10.5.1
reduces us to the case $k=1$.  Any other edge of length at
least $\xi$ is covered by Claim 1.  So we may assume that all
top edges have length less than $\xi$.  And then the result
follows easily from Inequalities 10.5.3 and 10.5.6.

Take $k=3$.
If there is an
edge of length at least $\xi$ lying around only one of the $v_i$,
then Inequality 10.5.1 reduces us to the case $k=2$.
If an edge of length at least $\xi$
lies around exactly two of the $v_i$, then it
is an edge of two of the quasi-regular tetrahedra.  These
quasi-regular tetrahedra give $2(0.55)\,\pt$, and the quasi-regular
tetrahedra around the third vertex $v_i$ give $0.55\,\pt$ more.
If a top edge of length at least $\xi$ lies around all three
vertices, then one of the endpoints of the edge
lies in $\{v_1,v_2,v_3\}$, so the result follows from Claim 1.
Finally, if all top edges have length at most $\xi$, we use
Inequalities 10.5.3, 10.5.6, 10.5.7.

Take $k=4$.  Suppose there
is a top edge $e$ of length at least $\xi$.  If $e$
lies around only one of the $v_i$, we
reduce to the case $k=3$.  If it lies around two of them, then
the two quasi-regular tetrahedra along this edge give $2(0.55)\,\pt$
and the quasi-regular tetrahedra around the other two vertices
$v_i$ give another $2(0.55)\,\pt$.  If both endpoints of $e$ are
among the vertices $v_i$, the result follows from Claim 2.  This
happens in particular if $e$ lies around four vertices.  If $e$
lies around only three vertices, one of its endpoints is one of
the vertices $v_i$, say $v_1$.  Assume $e$ is not around $v_2$.
If $v_2$ is not adjacent to $v_1$, then Claim 1 gives the
result.  So taking $v_1$ adjacent to $v_2$, we adapt Claim 1,
by using Inequalities 10.5.1--10.5.7, to show that the eight quasi-regular
tetrahedra around $v_1$ and $v_2$ give $4(0.55)\,\pt$.
Finally, if all top edges have length at most $\xi$, we use
Inequalities 10.5.3, 10.5.6, 10.5.7.
\qed
\enddemo

\bigskip

\head 6. Limitations on Types\endhead

Recall that a vertex of a planar map has type $(p,q)$ if it
is the vertex of exactly $p$ triangles and $q$ quadrilaterals.
This section restricts the possible types that appear
in a decomposition star.

Let $t_4$ denote the constant $0.1317\approx 2.37838774\,\pt$.
Proposition 4.1.11 asserts that every
quad cluster $R$ satisfies $\tau(R)\ge t_4$.

\proclaim{Lemma 6.1}  The following eight types $(p,q)$ are impossible:
(1)  $p\ge 8$,
(2) $p\ge 6$ and $q\ge 1$,
(3) $p \ge 5$ and $q\ge 2$,
(4) $p \ge 4$ and $q\ge 3$,
(5) $p \ge 2$ and $q\ge 4$,
(6) $p \ge 0$ and $q\ge 6$,
(7) $p \le 3$ and $q=0$,
(8) $p \le 1$ and $q=1$.
\endproclaim

\demo
{Proof}  By Proposition \refy{4.1.3} and Calculation \refy{10.1.3},
a lower bound on the dihedral
angle of $p$ simplices and $q$ quadrilaterals is
$0.8638p+1.153 q$.   If the type exists, this constant must
be at most $2\pi$.  One readily verifies in Cases 1--6
that $0.8638p+1.153q >2\pi$.  By Proposition \refy{4.3} and
Calculation \refy{10.1.2},
an upper bound on the dihedral angle of $p$ triangles and $q$
quadrilaterals is $1.874445 p + 3.247 q$.  In Cases 7 and 8 this
constant is less than $2\pi$.  \qed
\enddemo

\proclaim{Lemma 6.2}  If the type of any vertex of a decomposition star
is one of $(4,2)$, $(3,3)$, $(1,4)$, $(1,5)$, $(0,5)$, $(0,2)$,
$(7,0)$, then the decomposition star scores less than 8\,\pt.
\endproclaim

\demo{Proof}  According to Table 5.1, we have $\tlp(p,q)> \squander$,
for $(p,q) = (4,2)$, $(3,3)$, $(1,4)$, $(1,5)$, $(0,5)$, or $(0,2)$.
By Lemma 3.2, the result follows in these cases.
Now suppose that one of the vertices has type $(7,0)$.
By the results of Part I, which treats the case in which all standard
regions are triangles, we may assume that the star
has at least one quadrilateral.  We then
have $\tau(D^*)\ge\tlp(7,0) + t_4 >\squander$.  The result
follows.  \qed
\enddemo

In summary of the preceding two lemmas, we find that we may
restrict our attention to the following types of vertices.

$$\matrix
   (6,0)&      &       &       &       \\
   (5,0)&(5,1) &       &       &       \\
   (4,0)&(4,1) &       &       &       \\
        &(3,1) &(3,2)  &       &       \\
        &(2,1) &(2,2)  &(2,3)  &       \\
        &      &(1,2)  &(1,3)  &       \\
        &      &       &(0,3)  &(0,4)  \\
\endmatrix
$$

\head 7. Properties of Planar Maps\endhead

\proclaim{Proposition 7.1}
Suppose that
$\sigma(D^*)\ge 8\,\pt$.
The planar map $L$ of $D^*$ has the following properties
    (without loss of generality):
\endproclaim

1.  The graph $G(L)$ has no loops or multiple joins.

2.  Each face of $L$ is a triangle or quadrilateral.

3.  $L$ has at least eight triangular faces.

4.  $L$ has at most six quadrilateral faces, and at least one.

5.  Each vertex has one of the following types:
    $(6,0)$, $(5,0)$, $(4,0)$, $(5,1)$, $(4,1)$, $(3,1)$, $(2,1)$,
        $(3,2)$, $(2,2)$, $(1,2)$, $(2,3)$, $(1,3)$, $(0,3)$, and
        $(0,4)$.

6.  If $C$ is a $3$-circuit in $G(L)$, then it bounds a triangular
    face.

7.  If $C$ is a $4$-circuit in $G(L)$, then one of the following is
    true:

    \quad {\it (a)} $C$ bounds some quadrilateral region,

    \quad {\it (b)} $C$ bounds a pair of adjacent triangles,

    \quad {\it (c)} $C$ encloses one vertex, and it has type $(4,0)$ or $(2,1)$.

8.  $ t_4 q + \sum_{v\in V} (\tlp(p_v,q_v)- t_4q_v) \le \squander$,
    for any collection $V$ of vertices
    in $L$ such that no two vertices of $V$ lie on
    a common face.

\demo{Proof}  A loop would give a closed geodesic on the
unit sphere of length less than $\pi$.  A multiple join
would give nonantipodal conjugate points on the sphere.
Property 2 is
the restriction of step 3 of the Kepler conjecture.
Property 3 follows from I.9.1 and Part II.  Property 4
follows from $7 t_4>\squander$.
If there are no quadrilaterals,
the problem has been solved in Part I.  The restrictions
on types were obtained in Section 6.  Property 6 is
established in Part I.  A 4-circuit encloses at most one
vertex by I.4.2.  If it encloses none, it gives  a quad
cluster or two tetrahedra.  Otherwise, it encloses a
vertex of type $(4,0)$ or type $(2,1)$.  Property
8 comes from Section 5.
\qed\enddemo

\head 8. Combinatorics\endhead

In Steps III and IV of the Kepler conjecture, we need to generate
all planar maps satisfying various lists of conditions.
Here we describe a computer algorithm, which has been implemented
in {\it Java}.  We describe the algorithm in a way
that it can be used for Step IV as well.

We assume that the planar maps satisfy the following conditions.

1.  There are no loops or multiple joins.

2.  Each face is a polygon.

3.  The graph has between 3 and $N$ vertices, for some explicit $N$.

4.  The degree at each vertex is at most $6$.

5.  If $C$ is a $3$-circuit in $G(L)$, then $C$ bounds a triangular face.

6.  If $C$ is a $4$-circuit in $G(L)$, then one side of $C$ contains
    at most $1$ vertex.

7.  There is a constant $T\ge 0$ and constants $t_n\ge0$ and
$t(p,q,r)\ge 0$, such that
    $$\sum_I t_n + \sum_V\tlp(p_v,q_v,r_v) < T,$$
    where the first sum runs over {\it finished\/} faces, and
    the second sum runs over any {\it separated\/} set of vertices.
    (The term {\it finished\/}
     will be described below.  For planar maps
    produced as output of the algorithm, every face will be finished.
    But at intermediate stages of the algorithm some will be unfinished.)
    Also, a separated set of vertices means a collection of vertices
    with the properties that no two lie on a common face, and that
    every face
    at each of the vertices is finished.  Here $p_v$ is the number
    of triangles, $q_v$ is the number of quadrilaterals,
    $r_v$ is the number of other regions at $v$, and $n=n_v=p_v+q_v+r_v$.
    Also, $n$ is the number of sides of the finished face.

There are additional properties that it might be helpful to impose,
but the ones stated are sufficient for the description of the algorithm.
In our context, $r=r_v=0$, $t_3=0$,
$T=\squander$, $ t_4=0.1317$ and $t(p,q,0)=\tlp(p,q)-t_4 q$.

We produce all maps satisfying the conditions 1--7,
by extending
maps already satisfying these properties.  We assume that we have
a stack of maps $\{ L_j : j\in J\}$ satisfying the conditions  1--7
such
that each face of each of these maps is labeled finished or unfinished.
Each planar map in the stack will have at least one unfinished face.
In one iteration of the algorithm, we pop one of the maps $L$ from the
stack, modify it in various ways to produce new maps, output any of the
new maps that are finished (meaning all faces are finished), and push the
remaining ones back on the stack.  When all maps have been popped from
the stack, we are guaranteed to have produced all finished maps satisfying
properties (1--7). Here are the details of the algorithm.

1.    Let $L$ be a planar
map that has been popped from the stack.  Fix any unfinished face $F$ of $L$
and any edge $e$ of $F$.
Label the vertices of $F$ consecutively $1,\ldots,\ell$, with $1$ and
$\ell$ the endpoints of $e$.  For each $m=3,\ldots,N$,
let $A_m$ be the set of all $m$-tuples
$(a_i)\in {\Bbb Z}^m$
satisfying $(1=a_1\le a_2\le\cdots\le a_m=\ell)$,
with $a_{m-1}\ne a_m$.

2.  For each $a\in A_m$, we draw a new $m$-gon $F_a$
along the edge $e$ of  $F$ as follows.
We construct the vertices $v(1),\ldots,v(m)$
inductively.  If $i=1$ or $a_i\ne a_{i-1}$, then we set $v(i)=$ vertex $a_i$
of $F$.  But if $a_i=a_{i-1}$, we add a new vertex $v(i)$ to the planar map.
The new face is to be drawn along the edge
$e$ over the face $F$.
(For example, when $\ell = 5$, the faces $F_a$ corresponding to
$a = (1,1,3,4,4,5)$ and $a=(1,1,1,5)$ are shown in Diagram \refy{8.1}.)
As we run over all $m$ and $a\in A_m$, we run over all possibilities
for the finished face along the edge $e$ inside $F$.

\smallskip
\gram|2.2|\refy{8.1}|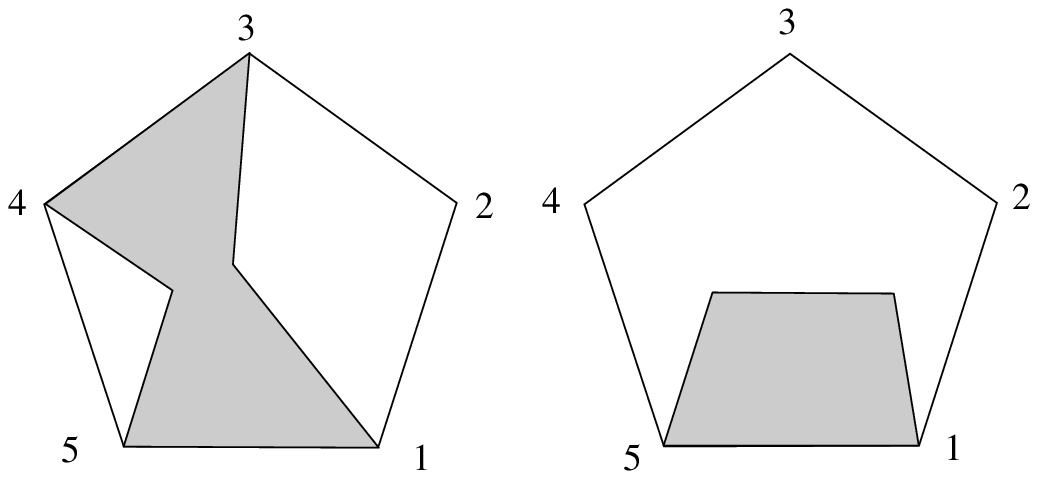|
\smallskip

3.  The face $F_a$ is to be marked as finished.  By drawing $F_a$, $F$
is broken into a number of smaller polygons.  (In Diagram \refy{8.1},
 $F$ is replaced respectively by three and two polygons.)
Each of these smaller polygons other than
$F_a$ is taken to be  unfinished, except triangles,
which are always taken to be
finished.

4.  Various planar maps extending $L$ are obtained by this process.
Those that do not satisfy the conditions (1--7) are discarded.
Those that have no unfinished faces are output.  The remaining ones
are pushed back onto the stack.  If the stack is empty, the
program terminates.  Otherwise, we pop a planar map from the
stack and return to the first step.   It is condition 7 that forces
the algorithm to terminate.

To begin the algorithm, we need an initial stack of planar maps.
The planar maps in the initial stack will be called {\it seeds}.
To produce a list of seeds, it is enough to give any list that is
guaranteed to produce all possibilities by the algorithm described
in $1,\ldots,4$.  For example, for Part III we want all configurations
with at least one quadrilateral and nothing but triangles and quadrilaterals.
We could let our initial stack consist of a single planar map $L$, where
the graph $G(L)$ is a single $4$-cycle.  Also, $L$ is to have two
faces, a finished quadrilateral and an unfinished complement.  Then
by iterating through the steps $1,\ldots,4$, we generate all possible
extensions of a quadrilateral to a planar map satisfying $1,\ldots,7$.

Although this seed would work,
in order to improve the performance of the algorithm,
in the implementation used for this paper,
we used a more detailed list of seeds, based on the classification
of types $(p,q)$ in Section 6.

This algorithm produced $1762$ cases,
even when the additional
properties listed in Proposition \refy{7.1} were used.
To be exact,  a few of the maps
may be superfluous, because there was no need to discard every
last map that we were allowed to.
The important point
is that an explicit finite list was obtained.  Because of the number
of possibilities involved we have not listed them here.  The
Java source code and pictures of the maps are available at \cite{H2}.

\head \refy{9}. Linear Programming Bounds\endhead

For each of the planar maps produced in Section \refy{8},
we define a linear programming
problem whose solution dominates the score of the decomposition stars
associated with the planar map.  A description of the linear
programs is presented in this section.

\proclaim{Theorem 9.1}  Let $L$ be any planar map obtained in
Section 8.  One of the following holds.  (1)  $L$ is the planar map
of the pentagonal prism, hexagonal-close packing, or face-centered
cubic.  (2) Every decomposition star with planar map $L$ scores
less than $8\,\pt$.  (3)  $L$ is one of the $18$ cases presented
in Appendix I.
\endproclaim

The variables of the linear
programming problem are the dihedral angles, the
scores of each of the standard clusters, and their edge lengths.

We subject
these variables to a system of linear inequalities.
First of all, the dihedral angles around each vertex sum to $2\pi$.
The dihedral angles, solid angles, and score are related by
the linear inequalities of Groups 1, 2, 3, and 4 in Section \refy{10}.
These include Propositions 4.1 and 4.2.
 The solid-angle variables
 are linear functions of dihedral angles.
The score of a decomposition star is
$$\sigma(S_1)+\cdots+\sigma(S_p)+\sigma(R_1)+\cdots+\sigma(R_q).$$
Forgetting the origin of the scores, solid angles, and dihedral
angles as nonlinear functions of the standard clusters and treating
them as formal
variables subject only to the given linear inequalities,
 we obtain a linear programming bound on the score.

Floating-point arithmetic was used freely in obtaining these
bounds. The linear programming package {\it CPLEX\/} was used (see
{\it www.cplex.com}). However, the results, once obtained, could
be checked rigorously as follows. (We did not actually do this
because the precision never seemed to be an issue, but this is how
it can easily be done. \footnote{The output from each linear
program in this paper has been double checked with interval
arithmetic. Predictably, the error bounds presented here were
satisfactory.  {\it 1/2002}})
 For each quasi-regular tetrahedron $S_i$
we have a nonnegative variable $x_i = \pt-\sigma(S_i)$.
For each quad cluster $R_k$, we have a nonnegative variable
$x_k = -\sigma(R_k)$.  A bound on the score is
$p\,\pt-\sum_{i\in I} x_i$, where $p$ is the number of triangular
standard regions, and $I$ indexes the faces of the planar map.
We give error bounds for a linear program involving scores and
dihedral angles.  Similar estimates can be made if there are
edges representing edge lengths.
Let the dihedral angles be $x_j$, for $j$ in some
indexing set $J$.  Write the linear constraints
as $Ax\le b$.  We wish to maximize $c\cdot x$ subject
to these constraints, where $c_i=-1$, for $i\in I$, and
$c_j=0$, for $j\in J$.  Let $z$ be an approximate solution
to the inequalities $zA\ge c$ and  $z\ge 0$ obtained by numerical
methods.  Replacing the negative entries of $z$ by $0$
we may assume that $z\ge0$ and that $zA_i> c_i-\epsilon$,
for $i\in I\cup J$, and some small error $\epsilon$.
If we obtain the numerical bound $p\,\pt+z\cdot b< 7.9999\,\pt$,
and if $\epsilon<10^{-8}$, then the score is less than $8\,\pt$.
In fact, note that
$$\left({z\over 1+\epsilon}\right) A_i$$
is at least $c_i$
for $i\in I$ (since $c_i=-1$),
and that it is greater than $c_i - \epsilon/(1+\epsilon)$,
for $i\in J$ (since $c_i=0$).
Thus, if $N\le 60$ is the number of vertices, and $p\le 2(N-2)\le116$
is the number of triangular faces,
$$
\align
\sigma(D^*) &\le p\,\pt + c\cdot x \le
         p\,\pt + \left({z\over 1+\epsilon}\right) A x
        + {\epsilon\over 1+\epsilon}\sum_{j\in J} x_j \\
    &\le p\,\pt + {z\cdot b\over 1+\epsilon} +
    {\epsilon\over 1+\epsilon} 2\pi N \\
    &\le \left[{p\,\pt+z\cdot b +
      {\epsilon}(p\,\pt+2\pi N)}\right]/(1+\epsilon)\\
    &\le \left[7.9999\,\pt +
        10^{-8}(116\,\pt+500)\right]/(1+10^{-8}) <8\,\pt.
\endalign
$$
\bigskip
In practice, we used $0.4429< 0.79984\,\pt$ as our cutoff, and
$N\le 14$ in the interesting cases, so much tighter error
estimates are possible.

\head \refy{10.} Calculations\endhead

In each of these calculations, when the cluster is a
quasi-regular tetrahedron $S$, we set
$\sigma=\sigma(S)$, $\dih=\dih(S)$, and so forth,
Let $\sigma_\lambda = \sigma-\lambda\zeta\pt\sol$,
for $\lambda = 1,3.2$.  We make similar abbreviations
for quad clusters.
The inequalities in Group 1 follow from
results appearing elsewhere.
 These inequalities have been
verified by interval arithmetic in \cite{H2}.

\define\n#1{\quad $#1.$\quad}

\subhead Group 1 \endsubhead
 Calculations that have been verified elsewhere.
{
\baselineskip = 0.66\baselineskip
\obeylines
\parskip=0pt

\hbox{}
{\it Quasi-regular tetrahedra: }
\n1  $\sigma\le\,\pt$   (I.9.1),
\n2  $\dih <  1.874445$ (I.8.3.2),
\n3  $\dih > 0.8638$   (I.9.3),
\n4  $\sigma <  -0.37642101\sol+0.287389$ (I.9.8),
\n5  $\sigma <  0.446634\sol-0.190249$ (I.9.9),
\n6  $\sigma <  -0.419351\sol+0.2856354+0.001$ (I.9.10,I.9.11,I.9.12,I.9.18).
\smallskip
{\it Quad clusters: }
\n7  $\sigma \le  0$   (II).
}

\bigskip
\subhead Group 2\endsubhead
Inequalities for quasi-regular tetrahedra depending
on edge lengths.
{
\baselineskip = 0.66\baselineskip
\obeylines
\parskip=0pt

\hbox{}
\n1  $\sol > 0.551285 + 0.199235(y_4+y_5+y_6-6)-0.377076(y_1+y_2+y_3-6)$,
\n2  $\sol < 0.551286 + 0.320937(y_4+y_5+y_6-6)-0.152679(y_1+y_2+y_3-6)$,
\n3  $\dih > 1.23095 -0.359894(y_2+y_3+y_5+y_6-8)+0.003(y_1-2)+0.685(y_4-2)$,
\n4  $\dih < 1.23096-0.153598(y_2+y_3+y_5+y_6-8)+0.498(y_1-2)+0.76446(y_4-2)$,
\n5  $\sigma <  0.0553737-0.10857(y_1+\cdots+y_6-12)$,
\n6  $\sigma+0.419351\sol <  0.28665-0.2(y_1+y_2+y_3-6)$,
\n7  $\sigma_1 <  10^{-6} -0.129119(y_4+y_5+y_6-6)-0.0845696(y_1+y_2+y_3-6)$.

}

\bigskip
\subhead Group 3\endsubhead
  General inequalities for quad clusters and quasi-regular
tetrahedra.
{
\baselineskip = 0.66\baselineskip
\obeylines
\parskip=0pt

\hbox{}
{\it Quasi-regular tetrahedra: }
\n1  $\sigma < 0.37898\dih -0.4111$,
\n2  $\sigma < -0.142\dih+ 0.23021$,
\n3  $\sigma < -0.3302\dih +0.5353$,
\n4  $\sigma_1 < 0.3897\dih -0.4666$,
\n5  $\sigma_1 < 0.2993\dih -0.3683$,
\n6  $\sigma_1 \le 0$,
\n7  $\sigma_1 < -0.1689\dih +0.208$,
\n8  $\sigma_1 < -0.2529\dih +0.3442$,
\n9  $\sigma_{3.2} < 0.4233\dih -0.5974$,
\n{10}  $\sigma_{3.2} < 0.1083\dih -0.255$,
\n{11}  $\sigma_{3.2} < -0.0953\dih -0.0045$,
\n{12}  $\sigma_{3.2} < -0.1966\dih +0.1369$,
\n{13}  $\sigma < -0.419351\sol + 0.796456\dih -0.5786316$,
\n{14}  $\sigma < -0.419351\sol + 0.0610397\dih +0.211419$,
\n{15}  $\sigma < -0.419351\sol  - 0.0162028\dih +0.308526$,
\n{16}  $\sigma < -0.419351\sol  - 0.0499559\dih +0.35641$,
\n{17}  $\sigma < -0.419351\sol  - 0.64713719\dih+ 1.3225$.
\smallskip
{\it Quad clusters: } Propositions 4.1 and 4.2.
}

\bigskip
\subhead Group 4\endsubhead
  Miscellaneous inequalities.
{
\baselineskip = 0.66\baselineskip
\obeylines
\parskip=0pt

\hbox{}
{\it Quasi-regular tetrahedra: }
\n1  The quasi-regular tetrahedra at a vertex of type $(4,0)$
    score at most $0.33\,\pt$ (I.5.2).
\n2  The sum of the dihedral angles around a vertex is $2\pi$.
\n3  The five quasi-regular tetrahedra $S_i$ at a vertex of type $(5,0)$
    satisfy $$\sum\sigma(S_i)  <  \sum (-0.419351\sol(S_i) + 0.2856354)$$
    (I.5.1.1).

\smallskip
{\it Flat quarters: }
\n4 $-0.398(y_2+y_3+y_5+y_6) + 0.3257y_1 - \dih_1 < -4.14938$,
    if $y_4\ge2.51$.
\n5 Proposition \refy{4.3}, Lemma \refy{5.3}.
\n6 Inequalities of Appendix 1: A.2.1--11, A.3.1--11, A.4.1--4, A.6.1--9,
    $A.6.1'$--$8'$, A.8.1--3.

}

\bigskip
\subhead Group 5\endsubhead
  Inequalities used by Lemma 5.3.
{
\baselineskip = 0.66\baselineskip
\obeylines
\parskip=0pt

\hbox{}
{\it quasi-regular tetrahedra: } Let $\xi=2.1773$, $m=0.2384$.
\n1  If $y_4\ge \xi$, then $\tau > 0.55\,\pt$,
\n2  If $y_4,y_5\ge\xi$, then $\tau > 2(0.55)\,\pt$,
\n3  If $y_4\le \xi$, then $\tau > -0.29349 + m\dih$,
\n4  If $y_4,y_6\le\xi$, $y_5\ge\xi$, then $\tau>-0.26303+m\dih$,
\n5  If $y_6\ge\xi$, $y_4,y_5\le\xi$, then $\tau>-0.5565+m(\dih_1+\dih_2)$,
\n6  If $y_4,y_5,y_6\le\xi$, then $\tau>-2(0.29349)+m(\dih_1+\dih_2)$,
\n7  If $y_4,y_5,y_6\le\xi$, then $\tau>-3(0.29349)+m(\dih_1+\dih_2+\dih_3)$.

\smallskip
Now set $\xi=2.177303$, $m =0.207045$.
\n8  If $y_4\ge \xi$, then $\sigma < (1-0.48)\,\pt$,
\n9  If $y_4,y_5\ge\xi$, then $\sigma < (1-2(0.48))\,\pt$,
\n{10}  If $y_4\le \xi$, then $\sigma < 0.31023815 - m\dih$,
\n{11}  If $y_4,y_6\le\xi$, $y_5\ge\xi$, then $\sigma<0.28365-m\dih$,
\n{12}  If $y_6\ge\xi$, $y_4,y_5\le\xi$, then $\sigma<0.53852-m(\dih_1+\dih_2)$,
\n{13}  If $y_4,y_5,y_6\le\xi$, then $\sigma<-pt+2(0.31023815)-m(\dih_1+\dih_2)$,
\n{14}  If $y_4,y_5,y_6\le\xi$, then $\sigma<-2\,\pt+3(0.31023815)-m(\dih_1+\dih_2+\dih_3)$.
}

\bigskip
\vfill\eject
\Refs
\bigskip

[F].  S.P. Ferguson, T.C. Hales, A Formulation of the Kepler
    Conjecture, preprint.

[I].  T.C. Hales, Sphere Packings I, Discrete and
    Computational Geometry, 17:1-51 (1997).

[II]. T.C. Hales, Sphere Packings II, Discrete and
    Computational Geometry, 18:135-149 (1997).

[V]. S.P. Ferguson, Sphere Packings V, thesis,
    University of Michigan, 1997.

[H1]. T.C. Hales, the Sphere Packing Problem, J. of Comp. and App. Math. 44
    (1992) 41--76.

[H2]. T.C. Hales, Packings
    {\tt http://www.math.lsa.umich.edu/\~\relax
        hales/packings.html}

[T]. W.T. Tutte, Graph theory, Addison-Wesley, 1984.

\endRefs
\newpage

\head Appendix 1. Some Final Cases\endhead

The graphs are numbered in the archive from $0$ to $1761$.
There are three cases that are treated elsewhere:
$\PM(4,1679)$ is the face-centered cubic, $\PM(4,1672)$ is
the pentahedral prism, and $\PM(4,1640)$ is the hexagonal
close packing.

The body of this paper eliminates all but a couple of dozen planar maps.
(The explicit list appears in the archive \cite{H2}.)
 This
appendix indicates how to eliminate the final cases.
The archive contains a few graphs that are isomorphic to each other.
The following discussion assumes that these duplicates have been
eliminated.
To exploit the nonlinearities of the problem, we use a
branch-and-bound method and divide the domains of the optimization
problem into several thousand smaller sets.   A linear
programming bound of $8\,\pt$ is obtained in each case.  This
appendix lists all of the inequalities that have been used, and
gives a description of the cases.  We refer the reader to \cite{H2} for
details about computer implementation of the linear programs.

\subhead A.1. Types\endsubhead
Each quad cluster with corners $(v_1,v_2,v_3,v_4)$ is one of four
types \cite{F}.
Although it was advantageous to group these cases together to simplify
the combinatorics, it is now better to separate these cases and to
develop linear inequalities for each case.

    1.  Two flat quarters with diagonal $(v_1,v_3)$.  The score of
    each quarter is compression or the analytic Voronoi function.

    2.
      Two flat quarters with diagonal $(v_2,v_4)$.  The score of
    each quarter is compression or the analytic Voronoi function.

    3.  Four upright quarters forming an octahedron.  The score of
    each upright quarter is compression or the averaged analytic
    Voronoi function $\octavor(Q)=(\vor(Q)+\vor(\hat Q))/2$.
    In other words, $\sigma(Q)=(\mu(Q)+\mu(\hat Q))/2$.  We
    are using notation from \cite{F}.

    4.  One of various mixed quad clusters.  The score is at most
    $\vor_0$, the truncated Voronoi function at radius $t_0=1.255$.

\smallskip
\subhead A.2. Flat quarters\endsubhead
For flat quarters, we have the following inequalities that were
established by interval arithmetic.   The edge $y_4$ is taken
to be the diagonal of the flat quarter. Here $\sigma=\mu$.

{
\baselineskip = 0.66\baselineskip
\obeylines
\parskip=0pt

\hbox{}

\n1  $-\dih_2+0.35 y_2 - 0.15 y_1 - 0.15 y_3 +0.7022 y_5 - 0.17 y_4 > -0.0123$,
\n2  $-\dih_3+0.35 y_3 - 0.15 y_1 - 0.15 y_2 +0.7022 y_6 - 0.17 y_4 > -0.0123$,
\n3  $\dih_2-0.13 y_2 + 0.631 y_1 + 0.31 y_3 -0.58 y_5+0.413 y_4+0.025 y_6 %
    > 2.63363$,
\n4  $\dih_3-0.13 y_3 +0.631 y_1+0.31 y_2 -0.58y_6+0.413 y_4 +0.025 y_5 %
    > 2.63363$,
\n5 $-\dih_1 +0.714y_1-0.221 y_2-0.221 y_3+0.92 y_4-0.221y_5-0.221 y_6 %
    > 0.3482$,
\n6 $\dih_1-0.315 y_1 +0.3972 y_2 +0.3972 y_3 - 0.715 y_4 +0.3972 y_5 %
    +0.3972 y_6 > 2.37095$,
\n7 $-\sol-0.187 y_1 -0.187 y_2 -0.187 y_3 +0.1185 y_4 + 0.479 y_5 %
    +0.479 y_6 > 0.437235$,
\n8 $\sol+0.488 y_1 + 0.488 y_2 + 0.488 y_3  - 0.334 y_5 %
    -0.334 y_6 > 2.244$,
\n9 $-\sigma -0.159 y_1 - 0.081 y_2 - 0.081 y_3 - 0.133 y_5 - 0.133 y_6 %
    > -1.17401$,
\n{10} $\sigma < -0.419351\sol + 0.1448 + 0.0436(y_5+y_6-4) %
    + 0.079431\dih$,
\n{11} $\sigma < 10^{-6} -0.197 (y_4+y_5+y_6-2\sqrt{2}-4)$.

}

\subhead A.3. Upright quarters\endsubhead
The following inequalities for upright quarters
have been established by interval arithmetic.
The first edge is taken to be the upright diagonal.

\bigskip
{
\baselineskip = 0.66\baselineskip
\obeylines
\parskip=0pt

\hbox{}

\n1 $\dih_1 - 0.636 y_1 + 0.462 y_2 + 0.462 y_3 - 0.82 y_4 + 0.462 y_5 %
    +0.462 y_6 > 1.82419$,
\n2 $-\dih_1 + 0.55 y_1 - 0.214 y_2 - 0.214 y_3 + 1.24 y_4 - 0.214 y_5 %
    -0.214 y_6 > 0.75281$,
\n3 $\dih_2 +0.4 y_1 -0.15 y_2 + 0.09 y_3 +0.631 y_4 -0.57 y_5 +0.23 y_6 %
    >2.5481$,
\n4 $-\dih_2-0.454 y_1 + 0.34 y_2 +0.154 y_3 -0.346 y_4 +0.805 y_5 %
    > -0.3429$,
\n5 $\dih_3 +0.4 y_1 -0.15 y_3 + 0.09 y_2 +0.631 y_4 -0.57 y_6 +0.23 y_5 %
    > 2.5481$,
\n6 $-\dih_3 -0.454 y_1 +0.34 y_3 +0.154 y_2 -0.346 y_4 +0.805 y_6 %
    > -0.3429$,
\n7 $\sol +0.065 y_2 + 0.065 y_3 + 0.061 y_4 -0.115 y_5 -0.115 y_6 %
    > 0.2618$,
\n8 $-\sol-0.293 y_1 -0.03 y_2 -0.03 y_3 + 0.12 y_4 +0.325 y_5 +0.325 y_6 %
    > 0.2514$,
\n9 $-\sigma-0.054 y_2 -0.054 y_3 - 0.083 y_4 - 0.054 y_5 -0.054 y_6 %
    > -0.59834$,
\n{10} $\sigma < -0.419351\sol + 0.079431\dih2 +0.06904 -0.0846(y_1-2.8)$,
\n{11} If $y_2,y_3\le 2.13$, then $\sigma < 0.07(y_1-2.51) %
    -0.133(y_2+y_3+y_5+y_6-8) - 0.135 (y_4-2)$.

}

\bigskip
\subhead A.4. Truncated quad clusters\endsubhead
Let $\phi(h,t) = (4 - 2\doct h t (h + t))/3$.
Set $t_0=1.255$ and $\phi_0=\phi(t_0,t_0)$.
In the truncated case $\vor_0$, \cite{F} gives
$$\vor_0 = \phi_0\sol + \sum A(y_i/2)\dih_i - 4\doct\sum_R\quo(R),$$
with
    $\phi_0 = \phi(t_0,t_0)$, and
$$A(h) = (1-h/t_0)(\phi(h,t_0)-\phi(t_0,t_0)).$$
Let $R$ be the Rogers simplex $R(y_1/2,\eta(y_1,y_2,y_6),t_0)$.
The function $\quo(R)$ is defined in \cite{F.3.3}.
We have $\quo(R)\ge0$.   Let $\vor_0^A$ denote the truncated Voronoi
score of half the quad cluster, divided into two simplices along
a diagonal, obtained by applying the formula
for $\vor_0$ to the simplex.
The following inequalities hold by interval arithmetic:

{
\baselineskip = 0.66\baselineskip
\obeylines
\parskip=0pt

\hbox{}

\n1 $\dih-0.372 y_1 +0.465 y_2 +0.465 y_3 + 0.465 y_5 +0.465 y_6 %
    > 4.885$,
\n2 $-\vor^A_0 - 0.06 y_2 -0.06 y_3 -0.185 y_5 -0.185 y_6 > -0.9978$,
    provided $\dih<2.12$, and $y_1,y_2,y_3\le 2.26$,
\n3 $-\vor^A_0+0.419351 \sol^A < 0.3072$,
    provided $\dih<2.12$, and $y_1,y_2,y_3\le 2.26$,
\n4 $\quo + 0.00758 y_1 + 0.0115 y_2 + 0.0115 y_6 > 0.06333$.

}

Also, $A\ge0$, $A'\le0$, and $A''\ge0$, for $h\in[1,t_0]$.
If $\dih\in[\dih_{\min},\dih_{\max}]$, and $h\in[1,h_{\max}]$,
for some constants $\dih_{\min}$, $\dih_{\max}$, and $1<h_{\max}\le t_0$,
then setting $\lambda = (A(h_{\max})-A(1))/(h_{\max}-1)$, we obtain
the additional elementary inequalities for $Ad := A(h)\dih$.

{
\baselineskip = 0.66\baselineskip
\obeylines
\parskip=0pt

\hbox{}
\n5 $Ad - A(1)\dih \le \lambda (h-1) \dih_{\min}$,
\n6 $Ad \le (A(1)+\lambda (h-1))\dih_{\max}$.

}

We use linear programming methods to determine bounds $h_{\max}$,
$\dih_{\min}$, $\dih_{\max}$.  If $A x\le b$ is the system of
inequalities used in the linear programs in the main body of the paper,
then we obtain an upper bound on a variable $y$ by solving the
linear program
$\max y$ subject to the constraints $A x\le b$, and the constraint
that the sum of the variables $\sigma$ (that is, the linear variables
corresponding to
the score) is at least $8\,\pt$.

\bigskip
\subhead A.5. Linear programs\endsubhead

Consider one of the remaining cases $\PM=\PM(4,n)$.
  Suppose $\PM$ has $r$ quadrilateral
faces.  We run $4^r$ linear programs, depending on which type
$1$--$4$ of quad cluster from A.1 each quadrilateral face represents.
In each case, we add the additional linear inequalities
from A.2, A.3, or A.4
as appropriate.  Note that a few of these inequalities are only
conditionally true, so that the inequality can be used only if
it is known that the condition holds.  All of the planar maps
have bounds under $8\,\pt$ by this method, except for
$\PM(4,71)$, $\PM(4,118)$, $\PM(4,126)$, and $\PM(4,178)$.

\bigskip
\subhead A.6. Quasi-regular tetrahedra\endsubhead

  The first nine inequalities assume
that $y_4+y_5+y_6\le 6.25$.  Of these, the last two
are established
only under the additional assumptions $y_1,y_2,y_3\le2.13$.
The next eight inequalities assume
that $y_4+y_5+y_6\ge 6.25$.  Of these, the last three
assume that $y_1,y_2,y_3\le2.13$.

{
\baselineskip = 0.66\baselineskip
\obeylines
\parskip=0pt

\hbox{}
\n1 $ \sol + 0.377076 y_1 + 0.377076 y_2 + 0.377076 y_3 - 0.221 y_4 %
    - 0.221 y_5 - 0.221 y_6  > 1.487741 $,
\n2 $ 0.221 y_4 + 0.221 y_5 + 0.221 y_6 - \sol  > 0.76822 $,
\n3 $ \dih + 0.34 y_2 + 0.34 y_3 - 0.689 y_4 + 0.27 y_5 + 0.27 y_6  %
    > 2.29295 $,
\n4 $ - \dih + 0.498 y_1 + 0.731 y_4 - 0.212 y_5 - 0.212 y_6  %
    > 0.37884 $,
\n5 $ - \sigma - 0.109 y_1 - 0.109 y_2 - 0.109 y_3 - 0.14135 y_4 - %
    0.14135 y_5 - 0.14135 y_6  > -1.5574737 $,
\n6 $ - \sigma - 0.419351 \sol - 0.2 y_1 - 0.2 y_2 - 0.2 y_3 - 0.048 %
    y_4 - 0.048 y_5 - 0.048 y_6  > -1.77465 $,
\n7 $ \tau - 0.0845696 y_1 - 0.0845696 y_2 - 0.0845696 y_3 - 0.163 %
    y_4 - 0.163 y_5 - 0.163 y_6  > -1.48542 $,
\n8 $\dih + 0.27 y_2 + 0.27 y_3 - 0.689 y_4 + 0.27 y_5 + 0.27 y_6 %
    > 2.01295 $,  if $y_1,y_2,y_3\le 2.13$.
\n9 $-\sigma - 0.14135 y_1 - 0.14135 y_2 - 0.14135 y_3 - 0.14135 y_4 %
    - 0.14135 y_5 - 0.14135 y_6 > -1.7515737 $, if $y_1,y_2,y_3\le2.13$.

}

\bigskip

{
\baselineskip = 0.66\baselineskip
\obeylines
\parskip=0pt

\hbox{}

\n{1'} $ y_4 +y_5 +y_6  > 6.25 $,
\n{2'} $ \sol + 0.378 y_1 + 0.378 y_2 + 0.378 y_3 - 0.1781 y_4 - %
     0.1781 y_5 - 0.1781 y_6 > 1.761445 $,
\n{3'} $ - \sol - 0.171 y_1 - 0.171 y_2 - 0.171 y_3 + 0.3405 y_4 + %
     0.3405 y_5 + 0.3405 y_6 > 0.489145 $,
\n{4'} $ - \sigma - 0.1208 y_1 - 0.1208 y_2 - %
     0.1208 y_3 - 0.0781 y_4 - 0.0781 y_5 - 0.0781 y_6 > -1.2436 $,
\n{5'} $ - \sigma - 0.419351 \sol - 0.2 y_1 - 0.2 y_2 - 0.2 y_3 + 0.0106 y_4 + %
     0.0106 y_5 + 0.0106 y_6 > - 1.40816 $,
\n{6'} $ \sol + 0.356 y_1 + 0.356 y_2 + 0.356 y_3 - 0.1781 y_4 - 0.1781 y_5 - %
     0.1781 y_6 > 1.629445 $, if $y_1,y_2,y_3\le2.13$,
\n{7'} $ - \sol - 0.254 y_1 - 0.254 y_2 - 0.254 y_3 + 0.3405 y_4 + 0.3405 y_5 + %
     0.3405 y_6 > - 0.008855 $, if $y_1,y_2,y_3\le2.13$,
\n{8'} $ - \sigma - 0.167 y_1 - 0.167 y_2 - 0.167 y_3 - 0.0781 y_4 - 0.0781 y_5 - %
     0.0781 y_6 > -1.51017 $, if $y_1,y_2,y_3\le2.13$.

}

\subhead A.7. Branch and Bound methods \endsubhead

For each case that failed the tests of A.5,
we use a branch-and-bound method
as follows.  We pick 10 quasi-regular
tetrahedra in the configuration.  We divide the domain into $2^{10}$
cases by imposing the constraint $y_4+y_5+y_6\le6.25$ or $y_4+y_5+y_6\ge6.25$
at each quasi-regular tetrahedron.  Depending on which constraint
is picked, we add (to the inequalities already present)
the first group 1--7 or the second group 1'--4'
of inequalities. Whenever we write of the branch-and-bound
inequalities associated with various faces, we mean these
two groups of inequalities.  When $r$ faces are used, there are
$2^{r}$ linear programs.
Before running these we determine the maximum dihedral angles
and edge lengths by separate linear programs.  If the conditions
of A.4.2, A.4.3 are met, we add these equations. Most of the $2^{10}$
linear programs
give bounds under $8\,\pt$. The only planar maps that still
give bounds over $8\,\pt$ are $\PM(4,126)$ and $\PM(4,178)$.

\bigskip
\subhead A.8. Final cases\endsubhead

The case $\PM(4,178)$ is not difficult.  The maximum height of
a vertex is less than 2.13, by linear programming bounds.  Add the
inequalities from A.3 and A.6 that assume this condition.  With
these additional inequalities, the linear programming bound is
less than $8\,\pt$.

In the case $\PM(4,126)$, again we find that the vertices have heights
less than $2.13$, and we add the relevant inequalities.  The score
is less than $8\,\pt$ unless the quad clusters fall into the
pure Voronoi case or
mixed cases.  (That is, there are no flat quarters and
no octahedra.)  In these remaining cases the diagonals of the
quad clusters are at least $2\sqrt2$.

Fix a vertex of type $(4,1)$ in the planar map.  There are two
such vertices to choose from, but there is an isomorphism carrying
one to the other.  Let $v_1,v_2,v_3,v_4$ be the corners of
one of the quad clusters, with $v_1$ the chosen vertex of type $(4,1)$.
The linear programming upper bound on the dihedral angle of
the quad cluster along $(0,v_1)$ is less than $1.694$.
The linear programming upper bound on $|v_2|+|v_1-v_2|+|v_1-v_4|+|v_4|$
is less than $8.709$.
The following inequality, established by interval computations,
implies that $|v_2-v_4|\in[2\sqrt2,2.93]$.

\n1   $\dih>1.694$, if $y_1,y_2,y_3\le2.13$, $y_2+y_3+y_5+y_6\le8.709$,
        and $y_4\ge2.93$.

This allows us to add the following two interval computations to
the set of linear programming inequalities.  Both hold for simplices
satisfying $y_1,y_2,y_3\le2.13$, $y_4\in[2\sqrt2,2.93]$, $y_5,y_6\in[2,2.51]$.
The quad cluster is broken in two simplices, with these inequalities
holding on each simplex.

\n2  $\dih_2+0.59 y_1+0.1y_2 +0.1y_3+0.55y_4 -0.6y_5-0.12y_6>2.6506$.

\n3  $\dih_2+0.35 y_1-0.24y_2+0.05y_3+0.35 y_4-0.72 y_5-0.18 y_6<0.47$.

With these additional inequalities, the linear programming bound
drops below $8\,\pt$.
We conclude that all $18$ planar maps that have scores below $8\,\pt$.
This concludes the proof of the main theorem of the paper and the
third step in the proof of the Kepler conjecture.

\newpage
\head Appendix 2. Interval Verifications\endhead

We make a few remarks in this appendix on the verification of
the inequalities of Proposition 4.1 and 4.2.
The basic method in proving an inequality $f(x)<0$ for $x\in C$,
is to use computer-based interval arithmetic to obtain rigorous
upper bounds on the second derivatives:
$|f_{ij}(x)|\le a_{ij}$, for $x\in C$.  These bounds lead immediately
to upper bounds on $f(x)$ through a Taylor approximation with
explicit bounds on the error.  We divide the domain $C$ as necessary
until the Taylor approximation
on each piece is less that the desired bound.

Some of the inequalities involve as many as 12 variables, such
as the octahedral cases of Proposition 4.2.  These are not directly
accessible by computer.  We describe some reductions we have used,
based on linear programming.
We start by applying the dimension reduction techniques described
in I.8.7.  We have used these whenever possible.

We will describe Proposition 4.2 because in various respects these
inequalities have been the most difficult to prove, although the
verifications of Propositions 4.1 and 4.3 are quite similar.
If there is a diagonal of length $\le2\sqrt{2}$, we have two flat
quarters $S_1$ and $S_2$.  The score breaks up into
    $\sigma=\sigma(S_1)+\sigma(S_2)$.  The simplices $S_1$ and
$S_2$ share a three-dimensional face.   The inequality we wish
to prove has the form
    $$\sigma \le a(\dih(S_1)+\dih_2(S_1)+\dih_2(S_2))+b.$$
We break the shared face into smaller domains on which we have
$$
\align
\sigma(S_1)&\le a (\dih(S_1)+\dih_2(S_1)) + b_1,\\
\sigma(S_2)&\le a \dih_2(S_1) + b_2,\\
\endalign
$$
for some $b_1,b_2$ satisfying $b_1+b_2\le b$.  These inequalities
are six-dimensional verifications.

If the quad cluster is an octahedron with upright diagonal, there
are four upright quarters $S_1,\ldots,S_4$.
We consider inequalities of the form
$$\sigma(S_i)\le \sum _{j\ne 4} a_j^i y_j(S_i)
    + a_7 (\dih_1(S_i)-\pi/2) + \sum_{j=2}^3 a \epsilon^i_j \dih_j(S_i)
    + b^i.\tag {$X_i$}$$
If $\sum_{i=1}^4 a_j^i \le0$, $j\ne 4$, and $\sum_i b^i\le b$, then
for appropriate $\epsilon^i_j\in\{0,1\}$, these inequalities
imply the full inequality for octahedral quad clusters.

By linear programming techniques, we were able to divide the
domain of all octahedra into about 1200 pieces and find inequalities
of this form on each piece, giving an explicit list of inequalities
that imply Proposition 4.2.  The inequalities involve six variables
and were verified by interval arithmetic.

To find the optimal
coefficients $a_j^i$ by linear programming we pose the linear problem
$$
\align
&\max t \\
&\hbox{such that}\\
&\quad X_i,\quad i=1,2,3,4, \ \{S_1,S_2,S_3,S_4\}\in C\\
&\quad\sum_i a_j^i \le 0,\\
&\quad\sum_i b^i \le b,
\endalign
$$
where $\{S_1,S_2,S_3,S_4\}$
runs over all octahedra in a given domain $C$.
The nonlinear inequalities $X_i$ are to be regarded as linear
conditions on the coefficients $a_j^i$, $b^i$, etc.
The nonlinear functions $\sigma(S_i),\dih(S_i)$, $y_j(S_i)$
are to be regarded as the coefficients of the variables $a_j^i,\ldots$
in a system of linear inequalities.  There
are
infinitely many constraints, because the set $C$ of octahedra is
infinite.  In practice we approximate
$C$ by a large finite set.  If the maximum of $t$ is positive,
then we have a collection of inequalities in small dimensions that
imply the inequality for octahedral quad clusters.  Otherwise,
we subdivide $C$ into smaller domains and try again.  Eventually,
we succeed in partitioning the problem into six-dimensional pieces,
which were verified by interval methods.

\smallskip
If the quad cluster is a mixed case, then \cite{F} gives
$$\sigma\le \vor_0, -1.04\,\pt,$$ so also
$$\sigma \le {3\over 4}\vor_0 + {1\over 4} (-1.04\,\pt).$$
In the pure Voronoi case with no quarters and no enclosed vertices,
we have the approximation
$$\sigma \le \vor(\cdot,\sqrt2) \le \vor_0.$$
If we prove $\vor_0\le a (\dih_1+\dih_2) + b$, the mixed case is established.
This is how Inequality 4.2.1 was established.  Dimension
reduction reduces this to a seven-dimensional verification.  We draw the
shorter of the two diagonals between corners of the quad cluster.
As we begin to subdivide this seven-dimensional domain, we are able to
separate the quad cluster into two simplices along the diagonal, each
scored by $\vor_0$.  This reduces the dimension further, to make
it accessible.  The two cases, 4.2.2 and 4.2.3, are similar, but we
establish the inequalities
$$
\align
{3\over 4}\vor_0 + {1\over 4} (-1.04\,\pt) &\le a (\dih_1+\dih_2) + b,\\
\vor(\cdot,\sqrt2)&\le a (\dih_1+\dih_2) + b.
\endalign
$$
This completes our sketch of how the verifications were made.

\bye